# Symplectic Geometry for Engineers – Fundamentals

Stefan Gössner

*Department of Mechanical Engineering, University of Applied Sciences, Dortmund, Germany.*

Dec. 2019

## 1. Introduction

The most simple standard example in Hamiltonian mechanics is the one-dimensional motion of a single particle. At a given time $t$ its state is determined by two real numbers – position $q$ and momentum $p$. The vector space $\mathbb{R}^2$ of positions and momenta as a trajectory is called *phase space*.

This is one famous example of several others, that lead to the guiding principle of modern physics – *physics is geometry* [1][2]. The mathematical theory underlying Hamiltonian mechanics is called *symplectic geometry* [2]. So symplectic geometry arose from the roots of mechanics and is seen as one of the most valuable links between physics and mathematics today [21].

> "*Symplectic geometry is an even dimensional geometry. It lives on even dimensional spaces, and measures the sizes of 2-dimensional objects rather than the 1-dimensional lengths and angles that are familiar from Euclidean and Riemannian geometry*" [3].
> (Dusa Mc. Duff)

Symplectic geometry in its simplest possible case is geometry of the plane $\mathbb{R}^2$. Surprise is, we can hardly find relevant publications regarding that simple 2D case. The reason for this might be, that higher even-dimensional spaces are by far more interesting and challenging for physicists and mathematicians.

For a first overview the excellent paper *Symplectization of Science* from Mark Gotay and James Isenberg [1] is highly recommended to the interested reader.

Engineers confronted with different real world problems in statics, kinematics, dynamics and electromagnetism realize, that they are frequently able to assign them to the field of two-dimensional problems. So from this point of view it seems worthwhile to start a more indepth discussion around that symplectic geometry thing. Primary goal is to estimate the suitability of symplectic geometry in $\mathbb{R}^2$ for solving standard mechanical problems in engineering education.

This paper is structured as follows. At the beginning we want to have a more abstract mathematical view at general $n$-dimensional vector spaces $\mathbb{R}^n$. This section can be easily skipped by less interested readers, who might go directly to section 3 then, where we will discuss the relevant planar case $\mathbb{R}^2$. In section 4 the *areal* character of planar symplectic geometry is highlighted with several very fundamental example problems in Euclidean geometry, kinematics and dynamics.

# 2. Symplectic Vector Spaces $\mathbb{R}^n$

We establish a general $n$-dimensional real vector space. A vector $\mathbf{u}$ is an ordered list of real numbers $u_1$ to $u_n$ then, which we will write in standard matrix notation as a column vector

$$\mathbf{u} = \begin{pmatrix} u_1 \\ \vdots \\ u_n \end{pmatrix}.$$

> **Definition 2.1 – Real vector space**
> Let $V = \mathbb{R}^n$ be a real vector space of dimension $n$.

The only thing we can do with elements of $V$ right now is adding them as well as scalar multiplying them by real numbers. We do this according to the usual rules of arithmetic.

Please note, that we cannot measure lengths and angles with our vector space yet. Being able to do so, we need to move over to an Euclidean space. Be patient.

> **Definition 2.2 – Bilinear form**
> A bilinear form on $V$ is a function $f: V \times V \to \mathbb{R}$ that takes each ordered pair of elements of $V$ to a real number. It is linear in each variable, when the other one is fixed and is
> - *non-degenerate*, if $f(\mathbf{u}, \mathbf{v}) = 0$ for all $\mathbf{v} \in V$ then $\mathbf{u} = \mathbf{0}$.
> - *symmetric*, if $f(\mathbf{u}, \mathbf{v}) = f(\mathbf{v}, \mathbf{u})$ for all $\mathbf{u}, \mathbf{v} \in V$.
> - *skew-symmetric*, if $f(\mathbf{u}, \mathbf{v}) = -f(\mathbf{v}, \mathbf{u})$ for all $\mathbf{u}, \mathbf{v} \in V$.
> - *positive definite*, if $f(\mathbf{u}, \mathbf{u}) > 0$ for all $\mathbf{u} \neq \mathbf{0}$.

Every bilinear form on $V = \mathbb{R}^n$ for all $\mathbf{u}, \mathbf{v} \in V$ has the shape

$$f(\mathbf{u}, \mathbf{v}) = \mathbf{u}^T \mathbf{A} \mathbf{v} \tag{2.1}$$

for some $n \times n$ matrix $\mathbf{A}$. The non-degeneracy here is equivalent to have a nonsingular matrix $\mathbf{A}$.

> **Definition 2.3 – Inner product space / Euclidean space**
> An inner product space is a vector space $V$ along with an inner product on $V$. An inner product on $V$ is a positive definite symmetric bilinear form $g$ on $V$. For $\mathbf{u}, \mathbf{v} \in V$ it is defined by
> $$g(\mathbf{u}, \mathbf{v}) = \mathbf{u}_1 \mathbf{v}_1 + \cdots + \mathbf{u}_n \mathbf{v}_n$$
> A real vector space with an inner product is also called *Euclidean space*.

An inner product is a generalization of – and for real vector spaces equivalent to – the dot product, which is denoted

$$\mathbf{u} \cdot \mathbf{v} = g(\mathbf{u}, \mathbf{v}). \tag{2.2}$$

Matrix $\mathbf{A}$ in relation (2.1) is the *identity matrix* $\mathbf{I}_n$ then. Using the inner product we are now able to define the *length of a vector* via its Euclidean norm.

$$\|\mathbf{u}\| = \sqrt{\mathbf{u} \cdot \mathbf{u}} \tag{2.3}$$

We also have a relationship between the dot product and the angle $\theta$ between two vectors by

$$\mathbf{u} \cdot \mathbf{v} = \|\mathbf{u}\|\|\mathbf{v}\| \cos \theta \,. \tag{2.4}$$

Please note that we cannot tell the angle's direction – from $\mathbf{u}$ to $\mathbf{v}$ or otherwise round – because of $\mathbf{u} \cdot \mathbf{v} = \mathbf{v} \cdot \mathbf{u}$. The special case of $\cos \theta = 0$ enables us to test two vectors for perpendicularity. So two vectors are orthogonal, if their dot product is vanishing.

> **Definition 2.4 – Symplectic vector space**
> A symplectic vector space is a vector space $V$ with a symplectic structure. A symplectic structure $\omega$ on $V$ is a non-degenerate skew-symmetric bilinear form on $V$. That $\omega$ is called *linear symplectic structure* and is defined for $\mathbf{u}, \mathbf{v} \in V$ by
>
> $$\omega(\mathbf{u}, \mathbf{v}) = -\omega(\mathbf{v}, \mathbf{u}) \tag{2.5}$$

> **Definition 2.5 – Compatible complex structure**
> A complex structure $\mathbf{J}$ on a symplectic vector space $(V, w)$ together with an inner product $g$
>
> $$\omega(\mathbf{u}, \mathbf{J}\mathbf{v}) = g(\mathbf{u}, \mathbf{v}) \tag{2.6}$$
>
> is called an $\omega$-compatible structure $J : V \to V$ such that $\mathbf{J}^2 = -\mathbf{I}_n$ [6].

An equivalent condition for $\mathbf{J}$ being compatible with $\omega$ is

$$\omega(\mathbf{J}\mathbf{u}, \mathbf{J}\mathbf{v}) = g(\mathbf{J}\mathbf{u}, \mathbf{v}) = g(\mathbf{v}, \mathbf{J}\mathbf{u}) = \omega(\mathbf{v}, \mathbf{J}^2\mathbf{u}) = \omega(\mathbf{v}, -\mathbf{u}) = \omega(\mathbf{u}, \mathbf{v}) \,. \tag{2.7}$$

Calling $\omega$ a *symplectic area*, equation (2.7) basically says that a transformation by $\mathbf{J}$ preserves symplectic area [3]. Since $\mathbf{J}$ is an *orthogonal operator*, it is also compatible with the inner product such that

$$g(\mathbf{J}\mathbf{u}, \mathbf{J}\mathbf{v}) = g(\mathbf{v}, \mathbf{u}) \,. \tag{2.8}$$

Furthermore $\mathbf{J}$ is a skew-symmetric matrix with $\mathbf{J}^T = \mathbf{J}^{-1} = -\mathbf{J}$. So we have for the determinant of $n \times n$ matrix $\mathbf{J}$

$$\det(\mathbf{J}) = \det(\mathbf{J}^T) = \det(\mathbf{J}^{-1}) = (-1)^n \det(\mathbf{J}) \,.$$

Suppose that $n$ is an odd integer, that expression results in $\det(\mathbf{J}) = -\det(\mathbf{J})$, which is valid only for a singular matrix $\mathbf{J}$ and so contradicting non-degeneracy of the compatible symplectic structure. Thus $n = dim(V)$ must be *even*, which is an important characteristic with symplectic vector spaces. As a consequence of this the shape of the unit symplectic matrix $\mathbf{J}$ is

$$\mathbf{J} = \begin{pmatrix} 0 & -\mathbf{I}_m \\ \mathbf{I}_m & 0 \end{pmatrix} \text{ with } m = n/2 \,. \tag{2.9}$$

Here relation (2.1) comes out as

$$\omega(\mathbf{u}, \mathbf{v}) = \mathbf{u}^T \mathbf{J}^T \mathbf{v} \,. \tag{2.10}$$

We finally have a symplectic vector space, capable of doing the following four things:

Table 1: Operations with the symplectic vector space

| Operation | Notation | Description |
|---|---|---|
| Vector addition | $\mathbf{u} + \mathbf{v}$ | Addition of two vectors. |
| Scalar multiplication | $s\mathbf{u}$ | Multiplication of a scalar quantity and a vector quantity. |
| Dot product | $\mathbf{u} \cdot \mathbf{v}$ | Multiplication of two vector quantities. |
| Orthogonal operation | $\mathbf{Ju}$ | Applying a complex structure to a vector. |

Please note as an important outcome of this section, that symplectic geometry lives on even dimensional vector spaces only.

## 3. Symplectic Vector Space $\mathbb{R}^2$

In two-dimensional real vector space $\mathbb{R}^2$ we have elements with ordered pairs of real numbers $a_x$ and $a_y$, which we will continue to write in standard matrix notation as column vectors

$$\mathbf{a} = \begin{pmatrix} a_x \\ a_y \end{pmatrix}.$$

We already equipped our vector space with an inner product, which we prefer to write as dot product (we even won't always write the dot).

$$\mathbf{ab} = a_x\, b_x + a_y\, b_y \tag{3.1}$$

The $w$-compatible complex structure $\mathbf{J}$ is a skew-symmetric $2 \times 2$ matrix here on $\mathbb{R}^2$.

$$\mathbf{J} = \begin{pmatrix} 0 & -1 \\ 1 & 0 \end{pmatrix} \tag{3.2}$$

Applying it to a planar vector performs an orthogonal transformation, i.e. a counterclockwise rotation by $\pi/2$ on it. We want to place a compact *tilde* '~' symbol over a vector variable to indicate the perpendicular vector to a given vector. Borrowing it from vector space $\mathbb{R}^3$ (see below) we call it *tilde-operator* – visually reflecting the skew-symmetry of its operator matrix.

$$\tilde{\mathbf{a}} = \mathbf{Ja} = \begin{pmatrix} -y \\ x \end{pmatrix} \tag{3.5}$$

Getting the perpendicular vector to a given vector in practice is easy. We just need to exchange both components and then negate the new first component.

Using the complex structure together with the inner product gets us to the *standard symplectic form* (2.5) on $\mathbb{R}^2$

$$\tilde{\mathbf{a}}\mathbf{b} = a_x b_y - a_y b_x. \tag{3.6}$$

It gives us the *area* of the parallelogram spanned by two vectors $\mathbf{a}$ and $\mathbf{b}$. And we get a *directed* or *oriented* area *from* $\mathbf{a}$ *to* $\mathbf{b}$ due to its inherent antisymmetry $\tilde{\mathbf{a}}\mathbf{b} = -\tilde{\mathbf{b}}\mathbf{a}$. So a symplectic structure on the plane is just an *area form* [3].

Area form (3.6) is named differently in literature. It is either called *perp dot product* [14], *skew product* [15] or *pseudo-sclar product* [16]. Arnold [5] calls it *skew-scalar product* and Lim [12] named it *symplectic inner product*.

We want to stick with *symplectic inner product*, as it is emphasizing the fact, that we do not have another product. It is rather the *standard inner product* together with a compatible complex structure (2.6).

## 3.1 Rules and Identities

Applying the tilde-operator to vectors in $\mathbb{R}^2$ follows the rules:

$$\begin{aligned}
\tilde{\tilde{\mathbf{a}}} &= -\mathbf{a} \\
\tilde{\mathbf{a}}\mathbf{a} &= 0 \\
\tilde{\mathbf{a}}\mathbf{b} &= -\mathbf{a}\tilde{\mathbf{b}} \quad Antisymmetry \\
\tilde{\mathbf{a}}\tilde{\mathbf{b}} &= \mathbf{a}\mathbf{b}
\end{aligned} \quad (3.7)$$

We can treat vector equations like algebraic equations. They can be added, scaled, squared and orthogonalized by applying the tilde-operator to them. Dot multiplying a vector equation by a vector results in a scalar equation. Multiplying by a vector again yields a vector equation in turn and so on alternating.

Using these rules, we can derive some identities:

Table 2: Identities with the symplectic inner product

| Identity | Expression | |
|---|---|---|
| Jacobi | $\tilde{\mathbf{a}}(\tilde{\mathbf{b}}\mathbf{c}) + \tilde{\mathbf{b}}(\tilde{\mathbf{c}}\mathbf{a}) + \tilde{\mathbf{c}}(\tilde{\mathbf{a}}\mathbf{b}) = \mathbf{0}$ | (I.1) |
| Grassmann | $\tilde{\mathbf{a}}(\tilde{\mathbf{b}}\mathbf{c}) + \mathbf{b}(\mathbf{c}\mathbf{a}) - \mathbf{c}(\mathbf{a}\mathbf{b}) = \mathbf{0}$ | (I.2) |
| Lagrange | $(\tilde{\mathbf{a}}\mathbf{b})^2 + (\mathbf{a}\mathbf{b})^2 = a^2 b^2$ | (I.3) |
| Grassmann | $\tilde{\mathbf{a}}(\tilde{\mathbf{b}}\mathbf{a}) + \mathbf{b}a^2 - \mathbf{a}(\mathbf{b}\mathbf{a}) = \mathbf{0}$ | (I.4) |
| Binet-Cauchy | $(\tilde{\mathbf{a}}\mathbf{b})(\tilde{\mathbf{c}}\mathbf{d}) = (\mathbf{a}\mathbf{c})(\mathbf{b}\mathbf{d}) - (\mathbf{a}\mathbf{d})(\mathbf{b}\mathbf{c})$ | (I.5) |

Identity (I.4) results from (I.2) by substituting $\mathbf{c}$ by $\mathbf{a}$. We want to proof and geometrically interpret some of these identities below.

## 3.2 Relationship to Vector Space $\mathbb{R}^3$

The cross product $\mathbf{a} \times \mathbf{b} = \tilde{\mathbf{a}}\mathbf{b}$ in vector space $\mathbb{R}^3$ is closely related to the symplectic inner product, as the skew-symmetric *cross product matrix* $\tilde{\mathbf{a}}$ build from spatial vector $\mathbf{a}$ multiplied by spatial vector $\mathbf{b}$ contains the symplectic inner product (3.6) in its $z$-coordinate [23]. This is, where we were borrowing the *tilde-operator* from.

The result of the *vector triple product* $\mathbf{a} \times (\mathbf{b} \times \mathbf{c})$ used with spatial vectors $\mathbf{a}, \mathbf{b}, \mathbf{c} \in \mathbb{R}^3$, while $z$-coordinates are set to zero, is identical to the result of the symplectic expression $\tilde{\mathbf{a}}(\tilde{\mathbf{b}}\mathbf{c})$. This is where Grassmann identity (I.2) – also named BAC-CAB identity – comes from [23].

The *scalar triple product* $\mathbf{a} \cdot (\mathbf{b} \times \mathbf{c})$ used with spatial vectors $\mathbf{a}, \mathbf{b}, \mathbf{c} \in \mathbb{R}^3$, while $z$-coordinates are set to zero, vanishes due to orthogonality. So there is no equivalent expression in $\mathbb{R}^2$.

## 3.3 Relationship to Complex Numbers

The symplectic vector space is equipped with a complex structure $\mathbf{J}$ according to (2.6). This also holds in general for tangent vectors, so for vector $\mathbf{a} = (x,y)^T$ there is [3]

$$\mathbf{J}\frac{\partial \mathbf{a}}{\partial x} = \frac{\partial \mathbf{a}}{\partial y} \quad and \quad \mathbf{J}\frac{\partial \mathbf{a}}{\partial y} = -\frac{\partial \mathbf{a}}{\partial x} \,.$$

Complex structure $\mathbf{J}$ plays the role of the imaginary unit $i$ (see $\mathbf{J}^2 = -\mathbf{I}$) performing a counterclockwise rotation by $\pi/2$. We can write down the expression

$$c\,\mathbf{a} + d\,\tilde{\mathbf{a}} = \begin{pmatrix} cx - dy \\ cy + dx \end{pmatrix}, \tag{3.8}$$

which is a linear combination of the orthogonal basis $\mathbf{a}$ and $\tilde{\mathbf{a}}$. And it is identical to the result of the complex multiplication

$$a\,(c + id) = (x + iy)(c + id) = (cx - dy) + i(cy + dx) \quad with \quad a \in \mathbb{C}\,.$$

Expression (3.8) is interpreted geometrically as a *similarity transformation*. In a same way

$$\mathbf{a}\,(c\,\mathbf{I} + d\,\mathbf{J})^{-1} = \frac{c\,\mathbf{a} - d\,\tilde{\mathbf{a}}}{c^2 + d^2} = \frac{1}{c^2 + d^2}\begin{pmatrix} cx + dy \\ cy - dx \end{pmatrix}, \tag{3.9}$$

conforms to the complex division

$$\frac{a}{c + id} = \frac{x + iy}{c + id} = \frac{(cx + dy) + i(cy - dx)}{c^2 + d^2} \quad with \quad a \in \mathbb{C}\,.$$

## 3.4 Areal Type of Geometry

The symplectic inner product (3.6) represents the directed area $\tilde{\mathbf{a}}\mathbf{b}$ of a parallelogram constructed by adjacent sides $\mathbf{a}$ and $\mathbf{b}$. We can even interpret the standard inner product (3.1) as a *directed parallelogram area* $\tilde{\mathbf{b}}\tilde{\mathbf{a}} = \mathbf{a}\mathbf{b}$ built by the sides $\mathbf{b}$ tand $\tilde{\mathbf{a}}$ according to Fig. 1.

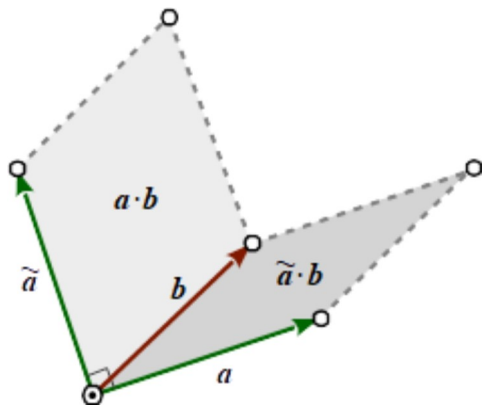

Fig. 1: Areas built by standard and symplectic inner product

So the triangles spanned by $\mathbf{a}$ and $\mathbf{b}$ as well as $\mathbf{b}$ and $\tilde{\mathbf{a}}$ are half of these areas each. *Directed areas* are *positive* if we go counterclockwise from $\mathbf{a}$ to $\mathbf{b}$ or from $\mathbf{b}$ to $\tilde{\mathbf{a}}$, otherwise *negative*.

## 3.5 Polar Vectors

Symplectic geometry being an *areal* type of geometry does not mean, that we cannot deal with *lengths* and *angles* as we are used to in Euclidean geometry. We might calculate the magnitude of a vector $\mathbf{a} = (x, y)^T$, using the Euclidean norm

$$a = \sqrt{\mathbf{a} \cdot \mathbf{a}} = \sqrt{x^2 + y^2}\,. \tag{3.10}$$

But we rather want to write it without the vertical bars. This way the magnitude can be written in a compact manner as the scalar component in its polar notation $\mathbf{a} = a\,\mathbf{e}_\alpha$. It does not even have to be positive.

There is no explicit division operation between two vectors defined. But an inverse element $\mathbf{a}^{-1}$ exists, so that $\mathbf{a}\mathbf{a}^{-1} = 1$, with

$$\mathbf{a}^{-1} = \frac{\mathbf{a}}{a^2}\,. \tag{3.11}$$

Now we are able to separate a vector $\mathbf{a}$ into its magnitude $a$ by using expression (3.10) and its unit direction vector $\mathbf{e}_\alpha$ and write down its polar notation

$$\mathbf{a} = a\,\mathbf{e}_\alpha = a\begin{pmatrix}\cos\alpha\\ \sin\alpha\end{pmatrix}. \tag{3.12}$$

Separating both vectors $\mathbf{a}$ and $\mathbf{b}$ into their polar components and reusing them in the standard and symplectic inner product, while applying trigonometric addition theorems, yields

$$\mathbf{a}\mathbf{b} = ab\,(\cos\alpha\cos\beta + \sin\alpha\sin\beta) = ab\cos(\beta - \alpha)$$

$$\tilde{\mathbf{a}}\mathbf{b} = ab\,(\cos\alpha\sin\beta - \sin\alpha\cos\beta) = ab\sin(\beta - \alpha)$$

Renaming the angular difference $\theta = \beta - \alpha$, we are able to get the *directed angle* from vector $\mathbf{a}$ to $\mathbf{b}$

$$\tan\theta = \frac{\tilde{\mathbf{a}}\mathbf{b}}{\mathbf{a}\mathbf{b}} \tag{3.13}$$

via the ratio of their areas in Fig 1. Squaring and adding both equations above leads us to Lagrange identity (I.3).

Please note that similarity transformation (3.8) in its special shape $\mathbf{a}\cos\varphi + \tilde{\mathbf{a}}\sin\varphi$ performs a pure rotation by angle $\varphi$ on vector $\mathbf{a}$.

## 4. Some Applications

We want to proof the usefulness of this vectorial approach with a few very basic geometric and mechanical examples.

## 4.1 Collinearity of Three Points

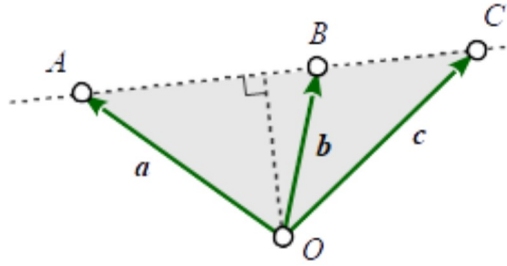

Fig. 2: Three collinear points

Let $\mathbf{a}, \mathbf{b}, \mathbf{c} \in \mathbb{R}^2$ be three vectors starting from a common point $O$. Testing their endpoints $A, B, C$ for collinearity, $(\tilde{\mathbf{b}} - \tilde{\mathbf{a}})(\mathbf{c} - \mathbf{b}) = 0$ must hold (Fig. 2). Outmultiplying results in

$$\tilde{\mathbf{a}}\mathbf{b} + \tilde{\mathbf{b}}\mathbf{c} + \tilde{\mathbf{c}}\mathbf{a} = 0, \qquad (4.1)$$

the condition of three vector endpoints to be collinear. This result is to be interpreted as the vanishing sum of three directed triangular areas. Observe the cyclic permutation pattern.

## 4.2 Simple Ratio

The segment length relation $\overline{AB} = \lambda \overline{BC}$ with respect to the three points $A, B, C$ in Fig. 2 is written as $(\mathbf{b} - \mathbf{a}) = \lambda(\mathbf{c} - \mathbf{b})$. Multiplying that equation by $\tilde{\mathbf{b}}$ results in the so called *simple ratio* [24]

$$\lambda = \frac{\tilde{\mathbf{a}}\mathbf{b}}{\tilde{\mathbf{b}}\mathbf{c}}. \qquad (4.2)$$

It is the ratio of two triangular areas. As both triangles have the same *height*, it is also the ratio of their baselines $\overline{AB}$ and $\overline{BC}$. If the baselines point in the same direction, the ratio is positive, otherwise negative. This conforms to the signs of the directed areas.

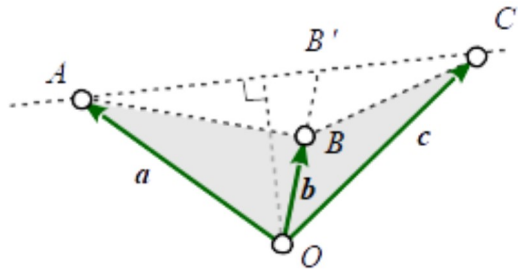

Fig. 3: Single Ratio

Simple ratio according to (4.2) even remains constant, when the endpoint $B$ of vector $\mathbf{b}$ isn't collinear with the other endpoints $A, C$ (Fig. 3). As a proof we scale vector $\mathbf{b}$ by $\mu$ to reach endpoint $B'$. So the condition for collinear vector endpoints $A, B', C$ now reads $(\mu\mathbf{b} - \mathbf{a}) = \lambda(\mathbf{c} - \mu\mathbf{b})$. Multiplying that expression by $\tilde{\mathbf{b}}$ results in expression (4.2) again, which is independent of factor $\mu$.

## 4.3 Double (Cross) Ratio

If four endpoints $A, B, C, D$ of vectors $\mathbf{a}, \mathbf{b}, \mathbf{c}, \mathbf{d}$ are collinear, a number is associated with them, which is called *double ratio* or *cross ratio* [24]. It is defined as the ratio of two simple ratios, i.e. ratio of four triangular areas.

$$\lambda = \overline{\frac{AC}{BC}} \Big/ \overline{\frac{AD}{BD}} = \frac{(\tilde{\mathbf{a}}\mathbf{c})(\tilde{\mathbf{b}}\mathbf{d})}{(\tilde{\mathbf{b}}\mathbf{c})(\tilde{\mathbf{a}}\mathbf{d})}.  \tag{4.3}$$

The double ratio defined by (4.3) plays an important role in projective geometry, where it is frequently noted as $[A, B; C, D]$.

## 4.5 Intersection of Two Lines

Point $A$ is element of line 1 with direction vector $\mathbf{u}$. Point $B$ lies on line 2 with direction vector $\mathbf{v}$. Let $\mathbf{a}$ be the vector from $A$ to $B$ (Fig. 3).

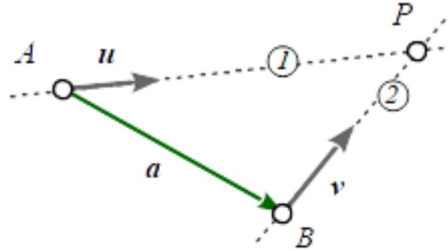

Fig. 4: Two intersecting lines

The loop closure equation of this triangle according to Fig. 4 then reads

$$\mathbf{a} + \mu\mathbf{v} - \lambda\mathbf{u} = \mathbf{0}.  \tag{4.4}$$

Multiplying this vector equation on the one hand by $\tilde{\mathbf{v}}$ and on the other hand by $\tilde{\mathbf{u}}$ helps to yield the two unknown scaling factors

$$\lambda = -\frac{\tilde{\mathbf{v}}\mathbf{a}}{\tilde{\mathbf{u}}\mathbf{v}} \quad and \quad \mu = \frac{\tilde{\mathbf{a}}\mathbf{u}}{\tilde{\mathbf{u}}\mathbf{v}},  \tag{4.5}$$

again as ratios of areas. With the knowledge of $\lambda$ we can now tell the vector from point $A$ to intersection point $P$ as $\lambda\mathbf{u}$. Introducing $\lambda$ and $\mu$ back into equation (4.4) and multiplying by common areal denominator $\tilde{\mathbf{u}}\mathbf{v}$ of expressions (4.5) leads to

$$(\tilde{\mathbf{u}}\mathbf{v})\mathbf{a} + (\tilde{\mathbf{v}}\mathbf{a})\mathbf{u} + (\tilde{\mathbf{a}}\mathbf{u})\mathbf{v} = \mathbf{0}.  \tag{4.6}$$

Please note the remarkable fact, that equation (4.6) as a triangle equation has the shape of the *Jacobi identity* (I.1). Thus Jacobi identity in $\mathbb{R}^2$ might be interpreted geometrically as a triangle equation.

In case of parallel lines denominators in expressions (4.5) are vanishing and intersection point $P$ goes to infinity. In the special case of a rectangular triangle, where $\mathbf{v} = \tilde{\mathbf{u}}$, point $P$ is the perpendicular projection of point $B$ onto line 1 in Fig. 4. Equation (4.6) will permutate into the shape of Grassmann identity (I.4) then.

## 4.6 Tangents to Two Circles

For two circles with radii $R_1$ and $R_2$ and distance vector $\mathbf{a}$ of their centers we want to determine their four tangents (Fig. 5).

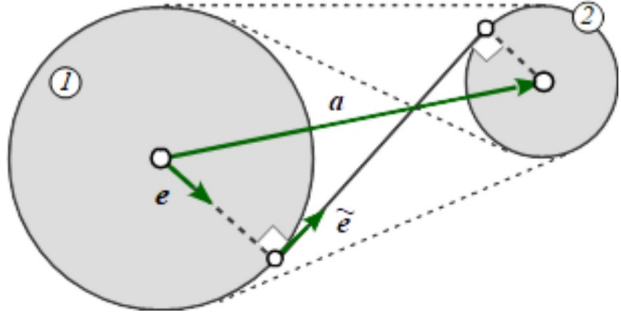

Fig. 5: Tangents to two circles

The loop closure equation reads

$$R_1 \mathbf{e} + \lambda \tilde{\mathbf{e}} \pm R_2 \mathbf{e} - \mathbf{a} = \mathbf{0}$$

Herein the positive sign of $R_2$ belongs to inner tangents and the negative sign to outer ones. Slightly reordering that equation and additionally writing down its 'orthogonal complement' results in

$$(R_1 \pm R_2)\mathbf{e} + \lambda \tilde{\mathbf{e}} = \mathbf{a}$$

$$(R_1 \pm R_2)\tilde{\mathbf{e}} - \lambda \mathbf{e} = \tilde{\mathbf{a}}$$

Squaring one of these equations leads us to the factor $\lambda$.

$$\lambda = \pm\sqrt{\mathbf{a}^2 - (R_1 \pm R_2)^2} \qquad (4.6)$$

Now multiplying the first equation by $(R_1 \pm R_2)$, the second one by $\lambda$ and subtracting the second from the first yields

$$[(R_1 \pm R_2)^2 + \lambda^2]\mathbf{e} = (R_1 \pm R_2)\mathbf{a} - \lambda\tilde{\mathbf{a}},$$

which we can resolve for unit direction vector $\mathbf{e}$

$$\mathbf{e} = \frac{(R_1 \pm R_2)\mathbf{a} - \lambda\tilde{\mathbf{a}}}{a^2} \qquad (4.7)$$

and we are done. The combination of two solutions each in (4.6) and (4.7) gives us the four tangents in question. If the expression under the root in (4.6) is negative, the circles are either overlapping (only two outer tangents) or one circle is completely inside of the other one (no tangents). Setting $R_2$ to zero gets us to the special case of two tangents from a point to a circle.

## 4.7 Kinematic Analysis

For solving kinematic problems we prefer to use the polar representation of vectors according to (3.10). Their magnitude and orientation are clearly separated then and may vary independently with the time, i.e.

$$\mathbf{r} = r\,\mathbf{e}_\varphi, \quad \dot{\mathbf{r}} = \dot{r}\,\mathbf{e}_\varphi + \dot{\varphi}\,r\,\tilde{\mathbf{e}}_\varphi, \quad \ddot{\mathbf{r}} = (\ddot{r} - \dot{\varphi}^2 r)\,\mathbf{e}_\varphi + (\ddot{\varphi}\,r + 2\dot{\varphi}\dot{r})\,\tilde{\mathbf{e}}_\varphi.$$

As an example we take an inverted slider crank mechanism (Fig.6).

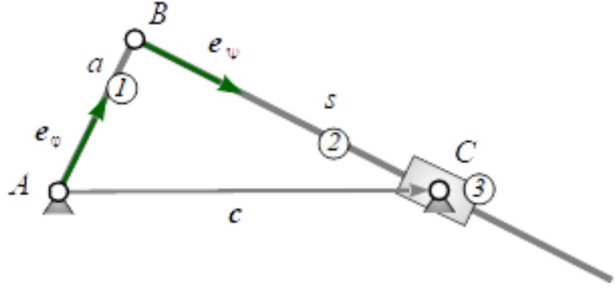

Fig. 6: Inverted slider crank

Crank 1 with fixed length $a$ is rotationally driven by $\dot\varphi = const$. Rocker 2 is connected to the crank by a rotational joint in point $B$ and guided by slider 3 fixed in point $C$ with variable orientation. We would like to know orientation $\psi$ and length $s$ as well as their first and second derivations with respect to time.

We start with loop closure equation (Fig. 6)

$$\mathbf{a} + s\,\mathbf{e}_\psi - \mathbf{c} = \mathbf{0}\,. \tag{4.8}$$

Resolving equation (4.8) for $s\mathbf{e}_\psi$ and squaring gets us to $s$ and allows to resolve for direction vector $\mathbf{e}_\psi$ then, i.e.

$$s = \sqrt{(\mathbf{c}-\mathbf{a})^2}\,,\quad \mathbf{e}_\psi = \frac{\mathbf{c}-\mathbf{a}}{s}\,. \tag{4.9}$$

Advancing to velocities we derive position equation w.r.t. time (4.8).

$$\dot\varphi\,\tilde{\mathbf{a}} + \dot s\,\mathbf{e}_\psi + \dot\psi\,s\,\tilde{\mathbf{e}}_\psi = \mathbf{0} \tag{4.10}$$

Multiplying (4.10) by unit vector $\mathbf{e}_\psi$ and by $\tilde{\mathbf{e}}_\psi$ respectively gives us the velocities of $s$ and $\psi$

$$\dot s = \dot\varphi\,\mathbf{a}\tilde{\mathbf{e}}_\psi\,,\quad \dot\psi = -\dot\varphi\,\frac{\mathbf{a}\mathbf{e}_\psi}{s}\,. \tag{4.11}$$

Please note both the standard and symplectic inner product in the velocity terms. Deriving velocity equation (4.10) with respect to time leads us to the accelerations. Crank 1 has no angular acceleration here.

$$-\dot\varphi^2\,\mathbf{a} + (\ddot s - \dot\psi^2 s)\,\mathbf{e}_\psi + (\ddot\psi + 2\dot\psi\dot s)\,\tilde{\mathbf{e}}_\psi = \mathbf{0} \tag{4.12}$$

Multiplying (4.12) also by unit vector $\mathbf{e}_\psi$ and by $\tilde{\mathbf{e}}_\psi$ respectively gives us the accelerations $\ddot s$ and $\ddot\psi$

$$\ddot s = \dot\psi(\dot\psi + \dot\varphi)s\,,\quad \ddot\psi = (\dot\varphi - 2\dot\psi)\frac{\dot s}{s}\,. \tag{4.13}$$

We were reusing velocity expressions (4.11) herein.

## 4.8 Hamiltonian Dynamics

The most simple example with Hamilton dynamics is the one-dimensional harmonic oscillator. A point mass $m$ is moving under the influence of the spring with ratio $k$ (Fig.7).

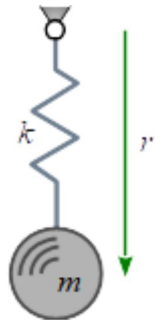

Fig. 7: Harmonic oscillator

We start with the Lagrangian $L = K - V$ as the difference of kinetic energy $K$ and potential energy $V$

$$L = \frac{1}{2}m\dot{r}^2 - \frac{1}{2}kr^2$$

and define the generalized coordinate $q$ and the generalized momentum $p$ by

$$q = r \quad \text{and} \quad p = \frac{\partial L}{\partial \dot{q}} = m\dot{q}. \tag{4.14}$$

It is an important point to accept these two variables as independent from each other. The pair of $q$ and $p$ uniquely describes the current state of the particle motion and by collecting them in a 2-dimensional vector $(q, p)^T$ we define a point in the phase space. Introducing $q$ and $p$ in the Lagrangian yields

$$L(q, p) = \frac{1}{2}\frac{p^2}{m} - \frac{1}{2}kq^2.$$

In order to advance from Lagrangian $L$ to the Hamiltonian $H$ we simply use the relation $H = p\dot{q} - L$, which is mathematically equivalent to applying the Legendre transformation.

$$H(q, p) = \frac{1}{2}\frac{p^2}{m} + \frac{1}{2}kq^2. \tag{4.15}$$

We recognize this expression as the total energy of the system, which is conservative, so the energy is constant. It can be shown in general that by variation of the Hamiltonian we get the relations

$$\dot{q} = \frac{\partial H}{\partial p} \quad \text{and} \quad \dot{p} = -\frac{\partial H}{\partial q}.$$

This can be encoded in terms of a symplectic transformation

$$\begin{pmatrix} \dot{q} \\ \dot{p} \end{pmatrix} = -\mathbf{J} \begin{pmatrix} \frac{\partial H}{\partial q} \\ \frac{\partial H}{\partial p} \end{pmatrix} = -\tilde{\nabla} H,$$

where $\nabla H$ is the gradient of $H$. The first vector component contains the velocity $\dot{q}$ and the second one Newton's equation of motion $m\ddot{q} = -kq$. Maurice de Gosson calls it "*just a fancy way to write Newton's second law*" [13].

We understand that the phase space as the trajectory built by all points $(q, p)^T$ – which can be shown to be an ellipse in this example – is even-dimensional also for systems with degrees of freedom greater than 1.

"The underlying idea of Hamiltonian mechanics is, that once one knows a particles initial state $(q_0, p_0)^T$ in phase space and also the forces that act on it, then one has enough information to chart the particle's motion in time" [1].

## 5. Conclusion

Symplectic geometry is a geometry of even dimensional spaces in which area measurements, rather than length measurements, are the fundamental quantities [7]. With symplectic geometry in $\mathbb{R}^2$ we have a *minimal apparatus* based on the four operations *addition*, *scalar multiplication*, *standard inner product* and *orthogonal operator*.

The standard inner product together with the orthogonal operator is the *symplectic inner product*. It has been shown, that both the symplectic as well as the standard inner product can be assigned a parallelogram area. So symplectic geometry in $\mathbb{R}^2$ is in fact a *pure areal type* of geometry. Additional benefits are:

- Symplectic geometry is coordinate free, i.e. coordinates aren't needed except we explicitly want them.
- The definition of an explicit origin is not required.
- Transformation matrices are rarely need. *Similarity transformation* mostly does what we want.
- In contrast to Euclidean geometry we can tell the direction of angles.

It was shown by a set of basic geometrical, kinematical and dynamical example problems, that symplectic geometry in $\mathbb{R}^2$ is suitable, beneficial and understandable for engineering education. The author has used this vectorial approach in the field of mechanism analysis, synthesis and design for several years now in research and student education. Results obtained with this method are advantageous and promising; they encourage to continue.

Please note that symplectic geometry is by far more than this minimal set of characteristics we were discussing here. For further in-depth information, reading of [5,6,10,11,13,17,19,20] is recommended.

## References

[1] M.J. Gotay and J.A. Isenberg, *The Symplectization of Science* [https://www.pims.math.ca/~gotay/Symplectization(E).pdf]
[2] D. McDuff (1999), *A Glimpse into Symplectic Geometry* [http://www.math.stonybrook.edu/~dusa/imu7.pdf]
[3] D. McDuff (2010), *What is Symplectic Geometry?* [http://barnard.edu/sites/default/files/ewmcambrevjn23.pdf]
[4] D. McDuff (1998), *Symplectic Structures — A New Approach to Geometry*, [http://www.math.stonybrook.edu/~dusa/noethmay22.pdf]
[5] V.I. Arnol'd, A.B. Givental', *Symplectic Geometry*, [https://www.maths.ed.ac.uk/~v1ranick/papers/arnogive.pdf]
[6] Z. Wang, Linear Symplectic Geometry [http://staff.ustc.edu.cn/~wangzuoq/Courses/15S-Symp/Notes/Lec01.pdf]
[7] M. Symington, *From Linear Algebra to the Non-squeezing Theorem of Symplectic Geometry*, [https://www.math.ias.edu/files/wam/SymingtonLecture1234.pdf]
[8] V. Guillemin, *Symplectic Techniques in Physics*, Cambridge University Press, 1990, ISBN: 978-0521389907
[9] R. Berndt, *Einführung in die Symplektische Geometrie*, Vieweg+Teubner, 1998, ISBN: 978-3528031022
[10] A. C. da Silva, *Lectures on Symplectic Geometry*, Lecture Notes in Mathematics, 1764, Springer-Verlag, Berlin, 2001. MR 1853077
[11] A. Weinstein, *Symplectic Geometry*, [https://projecteuclid.org/download/pdf_1/euclid.bams/1183548217]
[12] C. W. Lim, X. S. Xu, *Symplectic Elasticity: Theory and Applications*. Applied Mechanics Reviews, American Society of Mechanical Engineers, 2011, 63 (5), pp.050802.

[13] M. A. de Gosson, *The Symplectic Egg*. Applied Mechanics Reviews, American Society of Mechanical Engineers, 2011, 63 (5), pp.050802.

[14] Hill, F.S. Jr., *The pleasures of "perp dot" products*, Graphics Gems IV, pp.138-148. Academic Press, 1994.

[15] skew product [https://www.encyclopediaofmath.org/index.php/Skew_product]

[16] pseudo-scalar product [https://www.encyclopediaofmath.org/index.php/Pseudo-scalar_product]

[17] M. Zworski, *Basic Symplectic Geometry* [https://math.berkeley.edu/~zworski/symple.pdf]

[18] V.I. Arnold, *Symplectization Complexication and Mathematical Trinities*, [https://www.maths.ed.ac.uk/~v1ranick/papers/arnold4.pdf]

[19] P. Piccione, D.V. Tausk, *A Student's Guide to Symplectic Spaces, Grassmannians and Maslov Index*, [https://www.ime.usp.br/~piccione/Downloads/MaslovBook.pdf]

[20] K. Hartnett, *A Fight to Fix Geometry's Foundations*, Quanta Magazine, [https://www.quantamagazine.org/the-fight-to-fix-symplectic-geometry-20170209/]

[21] ThatsMath, Symplectic Geometry, [https://thatsmaths.com/2019/05/30/symplectic-geometry/]

[22] P. Hamill, *A Student's Guide to Lagrangians and Hamiltonians*, Cambridge University Press, 2013, ISBN: 978-1107617520

[23] Wikipedia, *Cross Product*, [https://en.wikipedia.org/wiki/Cross_product]

[24] Izu Vaisman, *Foundations of Three-Dimensional Euclidean Geometry*, Taylor & Francis Inc., 1980, ISBN: 978-0824769017